\documentclass{article}
\usepackage{amsmath,amsthm,amssymb}

\newtheorem{thm}{Theorem}
\newtheorem{prop}{Proposition}
\newtheorem{lem}{Lemma}
\newtheorem{defn}{Definition}

\newcommand{\qzeta}{\zeta_q}
\newcommand{\qzetabar}{\overline{\zeta}_q}
\newcommand{\C}{\mathbb{C}}
\newcommand{\Q}{\mathbb{Q}}

\newcommand{\depth}{\mathop\mathrm{dep}\nolimits}
\newcommand{\height}{\mathop\mathrm{ht}\nolimits}
\newcommand{\rank}{\mathop\mathrm{rank}\nolimits}
\newcommand{\Li}{\mathrm{Li}}
\newcommand{\bfk}{\boldsymbol{k}}

\let\epsilon\varepsilon

\newcommand{\be}{\begin{eqnarray*}} 
\newcommand{\en}{\end{eqnarray*}}
\newcommand{\bea}{\begin{eqnarray}} 
\newcommand{\ena}{\end{eqnarray}}
\newcommand{\nn}{\nonumber} 

\title{On relations for the $q$-multiple zeta values}
\author{OKUDA, Jun-ichi \and TAKEYAMA, Yoshihiro
\thanks{Research Fellow of the Japan Society for the Promotion of Science.}}

\date{}

\begin{document}
\maketitle

\begin{abstract} 
We prove some relations for the $q$-multiple zeta values ($q$MZV). 
They are $q$-analogues of 
the cyclic sum formula, the Ohno relation  
and the Ohno-Zagier relation for the multiple zeta values (MZV). 
We discuss the problem to determine 
the dimension of the space spanned by $q$MZV's over $\Q$, 
and present an application to MZV. 
\end{abstract}

\section{Introduction}
In this paper we prove some relations 
for a certain class of $q$-series 
called $q$-multiple zeta values ($q$MZV, for short). 

Let us recall the definition of $q$MZV \cite{math.QA/0304448}. 
We call a sequence of ordered positive integers $\bfk=(k_1,\dots,k_r)$ an index.
The weight, depth and height of the index are defined by
$|\bfk| := k_1+\dots+k_r$, $\depth \bfk := r$
and $\height \bfk := \#\,\{j|k_j\ge2\}$ respectively.
The index is said to be admissible if and only if $k_1\ge2$.

\begin{defn}
 For $0 < q < 1$ and 
 an admissible index $\bfk$,
 the $q$-multiple zeta value ($q$MZV) is defined by
 \begin{eqnarray}
 \qzeta(\bfk) :=
 \sum_{n_1>n_2>\dots>n_r>0}
  \frac{q^{n_1(k_1-1)+n_2(k_2-1)+\cdots+n_r(k_r-1)}}
       {[n_1]^{k_1}[n_2]^{k_2}\cdots [n_r]^{k_r}},
 \label{eq:def-qzeta} 
 \end{eqnarray}
 where $[n]$ is the $q$-integer 
 \[
  [n]:=\frac{1-q^n}{1-q^{\phantom{n}}}.
 \]
\end{defn}
Note that $0<1/[n]\le 1$ for any positive integer $n$. 
Hence the right hand side of \eqref{eq:def-qzeta} is absolutely convergent 
if $k_{1}>1$. 
In particular it is well defined as a $q$-series if $k_{1}>1$. 

By taking the limit $q \to 1$ of $q$MZV, we obtain 
the multiple zeta value (MZV, for short):
\[
 \lim_{q\uparrow1}\qzeta(\bfk) =
 \zeta(\bfk) :=
 \sum_{n_1>n_2>\dots>n_r>0}
  \frac{1}
       {n_1^{k_1}n_2^{k_2}\cdots n_r^{k_r}}.
\]
For MZV's there are many linear relations and algebraic ones 
over $\Q$.
The examples 
are the cyclic sum formula \cite{MR1971042}, the Ohno relation \cite{MR99k:11138} 
and the Ohno-Zagier relation \cite{MR2003e:11094}.
These relations do not suffice to give all relations of MZV's, 
and as a time the proof is much technical.
However these relations have quite beautiful structures
and it is valuable to construct explicit relations.

In the present paper, we 
give a $q$-analogue of the relations for MZV's and prove it.

First we show the $\Q$-linear relations for $q$MZV's
which are $q$-analogues of
the cyclic sum formula and
the Ohno relation:

\begin{thm}[The cyclic sum formula]\label{thm:cyclic}
 For any index $\bfk$ with some $k_i\ge2$,
 \begin{eqnarray*}
 \lefteqn{\sum_{i=1}^r
  \qzeta(k_i+1,k_{i+1},\dots,k_r,k_1,\dots,k_{i-1})}\\
&& =
 \sum_{i=1}^r
  \sum_{j=0}^{k_i-2}
   \qzeta(k_i-j,k_{i+1},\dots,k_r,k_1,\dots,k_{i-1},j+1).
 \end{eqnarray*}
\end{thm}

\begin{thm}[The Ohno relation]\label{thm:ohno}
 For any admissible index $\bfk=(k_1,\dots,k_r)$,
 there exit positive integers $a_1,b_1,a_2,b_2,\dots,a_s,b_s$
 such that
\begin{eqnarray}
  \bfk
  =
  (\underbrace{a_1+1,1,\dots,1}_{b_1},
    \underbrace{a_2+1,1,\dots,1}_{b_2},
     \dots,
      \underbrace{a_s+1,1,\dots,1}_{b_s}). \label{eq:def-code}
\end{eqnarray} 
Then the dual index $\bfk'=(k_1',\dots,k_{r'}')$ for $\bfk$ is defined by
\[
 \bfk'
  :=
  (\underbrace{b_s+1,1,\dots,1}_{a_s},
    \dots,
     \underbrace{b_2+1,1,\dots,1}_{a_2},
      \underbrace{b_1+1,1,\dots,1}_{a_1}).
  \nonumber
\]
For any admissible index $\bfk$ and non-negative integer $l$
\[
  \sum_{c_1+\dots+c_r=l \atop c_i\ge0}
   \qzeta(k_1+c_1,\dots,k_r+c_r)
  =
  \sum_{c'_1+\dots+c'_{r'}=l \atop c'_i\ge0}
   \qzeta(k'_1+c'_1,\dots,k'_{r'}+c'_{r'}).
\]
\end{thm}

These relations have the same form 
as the corresponding ones for MZV's. 
 
Secondly we show a $q$-analogue of the Ohno-Zagier relation: 

\begin{thm}[The Ohno-Zagier relation]\label{thm:ohno-zagier}
 We define a generating function of $q$MZV's as follows:
 \begin{eqnarray}
  \Phi_0(x,y,z)
  :=
  \sum_{k,r,s=0}^\infty
  \left\{
   \sum_{
     \bfk
     \atop
     |\bfk|=k,\ \depth\bfk = r,\ \height\bfk=s
    }
  \zeta_q(\bfk)
  \right\}
  x^{k-r-s} y^{r-s} z^{s-1}. \label{eq:generating-Ohno-Zagier}
 \end{eqnarray}
 Then
 \begin{eqnarray}
  \lefteqn{1+(z-xy)\Phi_{0}=}\nonumber\\
  &&\exp
   \left(
    \sum_{n=2}^{\infty}\qzeta(n) 
    \!
     \sum_{m=0}^{\infty}
      \frac{(q-1)^{m}}{m+n} (x^{m+n}+y^{m+n}-(\alpha^{m+n}+\beta^{m+n}))
   \right). \label{eq:q-Ohno-Zagier}
 \end{eqnarray}
 Here $\alpha^{m+n}+\beta^{m+n}$ is a polynomial in $x, y$ and $z$ 
 determined by
 \[
  \alpha+\beta=x+y+(q-1)(z-xy) \quad \hbox{and} \quad \alpha\beta=z. 
 \]
\end{thm}

By taking the limit $q \to 1$ we obtain the Ohno-Zagier relation 
for MZV's. 
The Ohno-Zagier relation for MZV's explains 
when the combinations of MZV's are in the algebra
generated by Riemann's zeta values $\zeta(k)$ over $\Q$.
Unfortunately
the product of $q$MZV's is not closed
in $\Q$-vector space spanned by $q$MZV's, 
however it is closed in $\Q[(1-q)]$-module
and preserves the weights
by counting the weight of $(1-q)$ by $1$ \cite{math.QA/0304448}.
{}From the consideration above we define 
the modified $q$MZV $\qzetabar(\bfk)$ by 
\begin{eqnarray} 
\qzetabar(\bfk):=(1-q)^{-|\bfk|}\qzeta(\bfk). 
\label{eq:modified-qMZV}
\end{eqnarray} 
Then the product of modified $q$MZV's is closed in 
the $\Q$-vector space spanned by them. 
Now the relation \eqref{eq:q-Ohno-Zagier} is rewritten as follows: 
\[
  1+(z-xy)\overline{\Phi}_{0}
  =
  \exp
   \left(
    \sum_{n=2}^{\infty}\qzetabar(n) 
     \sum_{m=0}^{\infty}
      \frac{(-1)^{m}}{m+n} (x^{m+n}+y^{m+n}-(\overline{\alpha}^{m+n}+\overline{\beta}^{m+n}))
   \right). 
\]
Here $\overline{\Phi}_{0}$ is defined by \eqref{eq:generating-Ohno-Zagier} 
with $\qzeta(\bfk)$ replaced by $\qzetabar(\bfk)$, and 
\[
 \overline{\alpha}+\overline{\beta}=x+y-(z-xy)
  \quad \hbox{and} \quad \overline{\alpha}\overline{\beta}=z. 
\]
Thus the relation \eqref{eq:q-Ohno-Zagier} explains when 
the linear combinations of the modified $q$MZV's are in the algebra 
generated by $\qzetabar(n) \,\, (n \ge 2)$ over $\Q$. 

In preparation of this paper we found \cite{math.QA/0402093}.
In \cite{math.QA/0402093},
Theorem~\ref{thm:cyclic}, Theorem~\ref{thm:ohno}
and Theorem~\ref{thm:ohno-zagier} of $z=0$ case
are proved independently.

The rest of the paper is organized as follows. 
In section \ref{sec:proofs} we prove the theorems above. 
Since the proofs are similar to the case of MZV, 
we omit some details presenting new features in the case of $q$MZV. 
In section \ref{sec:discussion} we discuss the problem to 
determine dimensions 
of certain spaces spanned by modified $q$MZV's over $\Q$. 
In the study of MZV it is an important problem 
to determine the dimension of the $\Q$-vector space spanned by 
MZV's of a fixed weight (see, e.g. \cite{MR96k:11110}). 
We consider a similar problem in the case of $q$MZV and 
present some observations. 

\section{Proofs}\label{sec:proofs}

\subsection{The Cyclic Sum Formula}

Set 
\begin{eqnarray*}
T(k_1,\dots,k_r) &:=&
 \sum_{n_1>\dots>n_r>n_{r+1}\ge0}
  \frac{q^{n_1-n_{r+1}} \cdot q^{n_1(k_1-1)+\dots+n_r(k_r-1)} }
       {[n_1-n_{r+1}][n_1]^{k_1}\cdots[n_r]^{k_r}},\\
S(k_1,\dots,k_r,k_{r+1}) &:=&
 \sum_{n_1>\dots>n_r>n_{r+1}>0}
  \frac{q^{n_1} \cdot q^{n_1(k_1-1)+\dots+n_r(k_r-1)+n_{r+1}(k_{r+1}-1)}}
       {[n_1-n_{r+1}][n_1]^{k_1}\cdots[n_r]^{k_r}[n_{r+1}]^{k_{r+1}}}.
\end{eqnarray*}
It is easy to see that 
$T(k_{1}, \ldots , k_{r})$ converges absolutely if 
all $k_{i}$'s are positive integers, and 
$S(k_{1}, \ldots , k_{r+1})$ converges absolutely if 
$\forall \, k_{i} \ge 1,$ or
$k_{r+1}=0$ and some of $k_i\ge2\ (i=1,\dots,r)$.

\begin{lem}\label{lem:cyclic}
Let $(k_{1}, \ldots , k_{r})$ be an index
with some $k_{i} \ge 2$. 
Then we have 
 \begin{eqnarray}
 \lefteqn{T(k_1,k_2,\dots,k_r) - \qzeta(k_1+1,k_2,\dots,k_r)}
  \nonumber\\
 &&=
 T(k_2,\dots,k_r,k_1)
 -\sum_{j=0}^{k_1-2} \qzeta(k_1-j,k_2,\dots,k_r,j+1).
  \label{eq:lem cyc}
 \end{eqnarray}
\end{lem}

We get Theorem~\ref{thm:cyclic}
by summing Lemma \ref{lem:cyclic} over all cyclic permutations of the 
sequence $(k_{1}, \ldots, k_{r})$.

\begin{proof}
 The left hand side of \eqref{eq:lem cyc} is
 \[
 T(k_1,\dots,k_r) - \qzeta(k_1+1,k_2,\dots,k_r)
 = S(k_1,\dots,k_r,0).
 \]
 Let us prove that $S(k_{1}, \ldots , k_{r}, 0)$ is equal to 
 the right hand side of \eqref{eq:lem cyc}. 

 We use the following partial fractional expansions:
 \begin{subequations}
 \begin{eqnarray}
  \frac{q^{n_1-n_{r+1}}}{[n_1-n_{r+1}][n_1]}
  &=&
  \left\{
  \frac{1}{[n_1-n_{r+1}]} - \frac{1}{[n_1]}
  \right\}
  \frac{1}{[n_{r+1}]}
  \label{eq:PFE}\\
  &=&
  \left\{
  \frac{q^{n_1-n_{r+1}}}{[n_1-n_{r+1}]} - \frac{q^{n_1}}{[n_1]}
  \right\}
  \frac{1}{[n_{r+1}]}. 
  \label{eq:PFE2}
 \end{eqnarray}
 \end{subequations}
 If $k_1\ge2$, by using \eqref{eq:PFE} we have
 \begin{eqnarray*}
 \lefteqn{S(k_1,k_2,\dots,k_r,k_{r+1})}\\
 &&=
 \sum_{n_1>\dots>n_r>n_{r+1}>0}
  \frac{q^{n_1}}{[n_1-n_{r+1}][n_1]}
  \frac{q^{n_1(k_1-1)+\dots+n_r(k_r-1)+n_{r+1}(k_{r+1}-1)}}
       {[n_1]^{k_1-1}[n_2]^{k_2}\cdots[n_r]^{k_r}[n_{r+1}]^{k_{r+1}}}\\
 &&=
 \sum_{n_1>\dots>n_r>n_{r+1}>0}
  \left\{
   \frac{1}{[n_1-n_{r+1}]} - \frac{1}{[n_1]}
  \right\}
   \frac{q^{n_1(k_1-1)+\dots+n_r(k_r-1)+n_{r+1}k_{r+1}}}
        {[n_1]^{k_1-1}[n_2]^{k_2}\cdots[n_r]^{k_r}[n_{r+1}]^{k_{r+1}+1}}\\
 &&=
 S(k_1-1,k_2,\dots,k_r,k_{r+1}+1) - \qzeta(k_1,k_2,\dots,k_r,k_{r+1}+1).
 \end{eqnarray*}
 By using this equality repeatedly we find 
 \begin{eqnarray}
 \lefteqn{S(k_1,k_2,\dots,k_r,0)}
 \nonumber\\
 &&
  =S(1,k_2,\dots,k_r,k_1-1)
   - \sum_{j=0}^{k_1-2} \qzeta(k_1-j,k_2,\dots,k_r,j+1). 
 \label{eq:pf-cyc1} 
 \end{eqnarray} 
 The equality above holds also in the case of $k_{1}=1$. 

 Now we consider the first term in the right hand side of \eqref{eq:pf-cyc1}:  
 \begin{eqnarray*}
 \lefteqn{S(1,k_2,\dots,k_r,k_{r+1})}\\
 &&=
  \sum_{n_1>\dots>n_r>n_{r+1}>0}
   \frac{q^{n_1}}{[n_1-n_{r+1}][n_1]}
   \frac{q^{n_2(k_2-1)+\dots+n_r(k_r-1)+n_{r+1}(k_{r+1}-1)}}
        {[n_2]^{k_2}\cdots[n_r]^{k_r}[n_{r+1}]^{k_{r+1}}}\\
 &&=
  \sum_{n_2>\dots>n_r>n_{r+1}>0}
   \frac{q^{n_2(k_2-1)+\dots+n_r(k_r-1)+n_{r+1}k_{r+1}}}
        {[n_2]^{k_2}\cdots[n_r]^{k_r}[n_{r+1}]^{k_{r+1}+1}}
  \sum_{n_1=n_2+1}^\infty
   \left\{
    \frac{q^{n_1-n_{r+1}}}{[n_1-n_{r+1}]} - \frac{q^{n_1}}{[n_1]}
   \right\}. 
 \end{eqnarray*}
 Here we used \eqref{eq:PFE2}. 
 Note that the sum $\sum_{n=1}^{\infty}q^{n}/[n]$ is convergent. 
 Hence we have 
 \[
 \sum_{n_1=n_2+1}^\infty
   \left\{
    \frac{q^{n_1-n_{r+1}}}{[n_1-n_{r+1}]} - \frac{q^{n_1}}{[n_1]}
   \right\}=
 \sum_{n_{r+1} > j > 0}\frac{q^{n_2-j}}{[n_2-j]}. 
 \]
 Therefore we find 
 \begin{eqnarray}
  S(1, k_{2}, \ldots , k_{r+1})=T(k_2,\dots,k_r,k_{r+1}+1).
 \label{eq:pf-cyc2} 
 \end{eqnarray}
 {}From \eqref{eq:pf-cyc1} and \eqref{eq:pf-cyc2} we obtain \eqref{eq:lem cyc}. 
\end{proof}

\subsection{The Ohno Relation}

The proof progresses as same as \cite{math.NT/0301277}.

For an admissible index $\bfk=(k_{1}, \ldots , k_{r})$, 
a sequence $(a_{1}, b_{1}, \ldots , a_{s}, b_{s})$ of positive integers 
is determined by the rule \eqref{eq:def-code}. 
We call the sequence a {\it code} of $\bfk$. 

Let $(a_1,b_1,\dots,a_s,b_s)$ be
the code of an admissible index $\bfk$.
We define generating functions of the left hand side and
the right hand side of the Ohno relation as follows:
\begin{eqnarray*}
 f(a_1,b_1,\dots,a_s,b_s;\lambda)
 &:=&
 \sum_{l=0}^\infty
  \sum_{c_1+\dots+c_r=l \atop c_i\ge0}
   \qzeta(k_1+c_1,\dots,k_r+c_r)
 \ \lambda^l,\\
 g(a_1,b_1,\dots,a_s,b_s;\lambda)
 &:=&
 \sum_{l=0}^\infty
  \sum_{c_1'+\dots+c_{r'}'=l \atop c_i'\ge0}
   \qzeta(k_1'+c_1',\dots,k_{r'}'+c_{r'}')
 \ \lambda^l\\
 &\hphantom{:}=& f(b_s,a_s,\dots,b_1,a_1;\lambda). 
\end{eqnarray*}
Here $\bfk'=(k_{1}', \ldots , k_{r'}')$ is the dual index of $\bfk$. 
These functions converge at $|\lambda|<1$,
and can be analytically continued to $\C\backslash\Omega$, where 
\[
 \Omega:=\{ q^{-n}[n] \, | \,  n \in \mathbb{Z}_{\ge 1} \},  
\]
by 
\begin{eqnarray}
 \lefteqn{f(a_1,b_1,\dots,a_s,b_s;\lambda)}
  \nonumber\\
  &&=
     \sum_{c_{1}, \ldots , c_{r}=0}^{\infty} 
    \sum_{n_1>\dots>n_r>0}
     \frac{q^{n_1(k_1-1)+\dots+n_r(k_r-1)}}
          {[n_1]^{k_1+c_1}\cdots [n_r]^{k_r+c_r}}\ \lambda^{c_{1}+\cdots +c_{r}}
 \nonumber\\
  &&=
  \sum_{n_1>\dots>n_r>0}
   \prod_{i=1}^s
     \frac{q^{n_{c_{i}}a_i}}{[n_{c_{i}}]^{a_i}}
   \prod_{j=1}^{r}\frac{1}{[n_{j}]-q^{j}\lambda},
  \label{eq:rewritten-f} 
\end{eqnarray}
 where $c_{1}=1$ and $c_i=b_1+\dots+b_{i-1}+1 \,\, (2 \le i \le s)$. 

Using these generating functions,
Theorem~\ref{thm:ohno} is stated as
\begin{eqnarray}
 f(a_1,b_1,\dots,a_s,b_s;\lambda)
 =
 g(a_1,b_1,\dots,a_s,b_s;\lambda)
 \label{eq:f=g}
\end{eqnarray}
for any code $(a_1,b_1,\dots,a_s,b_s)$.
To prove this equality, we prepare two propositions.

\begin{prop}\label{prop:analytic}
 \[
  f(a_1,b_1,\dots,a_s,b_s;\lambda)
  =
  \sum_{p=1}^\infty
  \frac{C_p}{[p]-q^p\lambda},
 \]
 where
 \[
  C_p
  =
  \sum_{d=1}^{r}
   \sum_{n_1>\dots>n_{d-1}>n_d=p \atop p>n_{d+1}>\dots>n_{r}}
      \prod_{i=1}^s 
      \frac{q^{n_{c_{i}}a_j}}
           {[n_{c_{i}}]^{a_j}}
     \prod_{j=1 \atop j\not=d}^{r}
     \frac{1}{[n_j]-q^{n_j-p}[p]}
 \]
\end{prop}

\begin{proof}
 {}From \eqref{eq:rewritten-f} we have 
 \begin{eqnarray*}
  \lefteqn{f(a_1,b_1,\dots,a_s,b_s;\lambda)}\\
  &&=
   \sum_{n_1>\dots>n_r>0}
    \prod_{i=1}^s
     \frac{q^{n_{c_{i}}{a_j}}}{[n_{c_{i}}]^{a_j}}
      \sum_{d=1}^{r} 
      \prod_{j\not=d}\left(\frac{1}{[n_j]-q^{n_j-n_i}[n_i]}\right) 
      \frac{1}{[n_d]-q^{n_d}\lambda}.
 \end{eqnarray*}

 If $\lambda$ is in a compact set in $\C\backslash\Omega$, we have 
 \begin{eqnarray*}
  \lefteqn{
   \sum_{n_1>\dots>n_r>0}\sum_{d=1}^{r}
    \left| 
      \prod_{i=1}^s
      \frac{q^{n_{c_{i}}{a_i}}}{[n_{c_{i}}]^{a_i}}
      \prod_{j\not=d}
      \left(
       \frac{1}{[n_j]-q^{n_j-n_i}[n_i]}
      \right) 
      \frac{1}{[n_d]-q^{n_d}\lambda}
    \right|}\\
  &&\qquad\qquad\le
  C
   \sum_{n_1>\dots>n_r>0}
    \prod_{i=1}^s
     \frac{q^{n_{c_{i}}{a_i}}}{[n_{c_{i}}]^{a_i}}
      \prod_{j=1}^{r}\frac{1}{[n_j]}\\
  &&\qquad\qquad=
  C\ \qzeta(k_1,\dots,k_r)<+\infty
 \end{eqnarray*}
 for a positive constant $C$. 
 Hence we can change the order of the summation
 and obtain the lemma.
\end{proof}

Before we state the second proposition (Proposition \ref{prop:difference} below) 
we introduce some conventions.  
We allow $0$'s appear in the code $(c_{1}, \ldots , c_{2s})$ 
with the identification 
\[
 (\dots,c_{i-1},0,c_{i+1},\dots) = (\dots,c_{i-1}+c_{i+1},\dots)
\]
for $2 \le i \le 2s-1$. 
It is consistent with the rule \eqref{eq:def-code}. 
By definition we set 
\[
 f(a_{1}, b_{1}, \ldots , a_{s}, b_{s} ; \lambda)=0 \quad 
 \hbox{if} \quad a_{1}=0 \,\, \hbox{or} \,\, b_{s}=0. 
\]

\begin{prop}\label{prop:difference}
Set $I:=\{(0,0), (0,1), (1,0), (1,1)\}$. 
For $\eta=\{(\epsilon_i, \delta_{i})\}_{i=1}^{s} \in I^{s}$ we set 
\[
 |\eta|=\sum_{i=1}^{s}(\epsilon_{i}+\delta_{i}), \quad 
 h(\eta)=\# \{ i \,\, | (\epsilon_{i}, \delta_{i})=(1, 1) \}. 
\]
Then we have 
 \begin{eqnarray}
  \lefteqn{
    \sum_{\eta \in I^{s}}
      (1-q)^{h(\eta)}(-\lambda)^{s-|\eta|+h(\eta)}
      f(a_1-\epsilon_1,b_1-\delta_1,\dots,a_s-\epsilon_s,b_s-\delta_s;\lambda)}
\nonumber\\
   && {}=
    \sum_{\eta' \in I^{s-1}}\sum_{\epsilon_{1}', \delta_{s+1}'=0, 1}
     (1-q)^{h(\eta')}
      (-\lambda')^{s-|\eta|-\epsilon_{1}'-\delta_{s+1}'+h(\eta)}
  \label{eq:difference}\\ 
   && \qquad {}\times  
       f(a_1-\epsilon_{1}', b_1-\delta_{2}',a_2-\epsilon_{2}',
        \dots,b_{s-1}-\delta_{s}',a_s-\epsilon_{s}',b_s-\delta_{s+1}';\lambda'),
	\nonumber
 \end{eqnarray}
 where $\lambda':= q\lambda-1$.
 In the left hand side $\eta=\{(\epsilon_{i}, \delta_{i})\}_{i=1}^{s}$, 
 and in the right hand side $\eta'=\{(\delta_{i}', \epsilon_{i}')\}_{i=2}^{s}$. 
\end{prop}

To prove this proposition, we use the following function:
\begin{eqnarray}
  \lefteqn{\rho((a_1,d_1),b_1,\dots,(a_s,d_s),b_s;\lambda)}
   \nonumber\\
  &&:=
  \sum_{n_1>\dots>n_{r}>0}
   \prod_{i=1}^s
    \frac{q^{(n_{c_{i}}-d_j)a_i}}
         {[n_{c_{i}}-d_j]^{a_i}}
   \prod_{j=1}^{r}\frac{1}{[n_{j}]-q^{n_{j}}\lambda}, 
\label{eq:def-rho}
\end{eqnarray}
where $r=\sum_{i=1}^{s}(a_{i}+b_{i})$, 
$c_{1}=1$ and $c_{i}=b_{1}+\cdots +b_{i-1}+1 \,\, (2 \le i \le s)$. 
Then the generating function $f$ is given by 
\[
 f(a_1,b_1,\dots,a_s,b_s;\lambda)=
 \rho((a_1,0),b_1,\dots,(a_s,0),b_s;\lambda).
\]

In the following we use the identification 
\begin{eqnarray*}
 (\dots,b_{i-1},(0,d_i),b_i,\dots) &=& (\dots,b_{i-1}+b_i,\dots),\\
 (\dots,(a_{i-1},d),0,(a_i,d),\dots) &=& (\dots,(a_{i-1}+a_i,d),\dots).
\end{eqnarray*}
It is consistent with the definition of $\rho$ \eqref{eq:def-rho}. 

\begin{lem}\label{lem:difference}
 The function $\rho$ satisfies following relations:

 \setlength{\leftmargini}{0cm}
 \setlength{\leftmarginii}{1.5em} 
 \begin{enumerate}
  \item
   \begin{enumerate}
   \item\label{item:1a}
    If $a_1\ge2$,
    \begin{eqnarray*}
      &&
      \lambda \rho((a_1,0),b_1,(a_2,0),\dots;\lambda)\\
      &&\quad
      - \rho((a_1-1,0),b_1,(a_2,0),\dots;\lambda)
       - \rho((a_1,0),b_1-1,(a_2,0),\dots;\lambda)\\
      &&\qquad
      -(1-q)\rho((a_1-1,0),b_1-1,(a_2,0),\dots;\lambda)\\
     &&=
     \lambda' \rho((a_1,1),b_1,(a_2,0),\dots;\lambda)
      - \rho((a_1-1,1),b_1,(a_2,0),\dots;\lambda).
    \end{eqnarray*}
   \item\label{item:1b}
    If $a_1=1$,
    \begin{eqnarray*}
      \lambda \rho((1,0),b_1,\dots;\lambda)
      - \rho((1,0),b_1-1,\dots;\lambda)
      =
      \lambda' \rho((1,1),b_1,\dots;\lambda).
    \end{eqnarray*}
   \end{enumerate}  
  \item \label{item:2}
  For $2\ge i \ge s-1$ or $i=s$ with $b_s\ge2$,
  \begin{eqnarray*}
    &\begin{array}{r@{\,}l}
    \lambda \rho(\dots,(a_{i-1},1),b_{i-1},&(a_i\hphantom{{}-1},0),b_i\hphantom{{}-1},(a_{i+1},0),\dots;\lambda)\\
    - \rho(\dots,(a_{i-1},1),b_{i-1},&(a_i-1,0),b_i\hphantom{{}-1},(a_{i+1},0),\dots;\lambda)\\
    - \rho(\dots,(a_{i-1},1),b_{i-1},&(a_i\hphantom{{}-1},0),b_i-1,(a_{i+1},0),\dots;\lambda)\\
      - (1-q)\rho(\dots,(a_{i-1},1),b_{i-1},&(a_i-1,0),b_i-1,(a_{i+1},0),\dots;\lambda)
   \end{array}\\
   &
   \begin{array}[t]{rr@{\,}l}
   =&\lambda' \rho(\dots,(a_{i-1},1),&b_{i-1}\hphantom{{}-1},(a_i\hphantom{{}-1},1),b_i,(a_{i+1},0),\dots;\lambda)\\
    &- \rho(\dots,(a_{i-1},1),&b_{i-1}\hphantom{{}-1},(a_i-1,1),b_i,(a_{i+1},0),\dots;\lambda)\\
    &- \rho(\dots,(a_{i-1},1),&b_{i-1}-1,(a_i\hphantom{{}-1},1),b_i,(a_{i+1},0),\dots;\lambda)\\
    &- (1-q)\rho(\dots,(a_{i-1},1),&b_{i-1}-1,(a_i-1,1),b_i,(a_{i+1},0),\dots;\lambda).
   \end{array}  
  \end{eqnarray*}
  \item
  \begin{enumerate}
   \item\label{item:3a}
    If $b_s\ge2$,
    \begin{eqnarray*}
      \lefteqn{\rho((a_1,1),b_1,\dots,(a_s,1),b_s;\lambda)}\\
      &&= \rho((a_1,0),b_1,\dots,(a_s,0),b_s;\lambda')
       -\frac{1}{\lambda'}\rho((a_1,0),b_1,\dots,(a_s,0),b_s-1;\lambda').
    \end{eqnarray*}
   \item\label{item:3b}
   If $b_s=1$,
   \begin{eqnarray*}
    &&
     \lambda\rho((a_1,1),b_1,\dots,(a_{s-1},1),b_{s-1},(a_s,0),1;\lambda)\\
    &&\quad
     -\rho((a_1,1),b_1,\dots,(a_{s-1},1),b_{s-1},(a_s-1,0),1;\lambda)\\
    &&=
    \lambda'\rho((a_1,0),b_1,\dots,b_{s-1},(a_s,0),1;\lambda')\\
    &&\quad
     -\rho((a_1,0),b_1,\dots,b_{s-1},(a_s-1,0),1;\lambda')\\
    &&\quad\quad
     -\rho((a_1,0),b_1,\dots,b_{s-1}-1,(a_s,0),1;\lambda')\\
    &&\quad\quad\quad
     -(1-q)\rho((a_1,0),b_1,\dots,b_{s-1}-1,(a_s-1,0),1;\lambda').
   \end{eqnarray*}
  \end{enumerate}
 \end{enumerate}
\end{lem}

\begin{proof}
Here we prove (\ref{item:1a}). 
The proofs of the other relations are similar. 

We have 
\begin{eqnarray}
 && \lambda \rho((a_1,0),b_1,\dots;\lambda)
   - \rho((a_1-1,0),b_1,\dots;\lambda)
  =
  \label{eq:pf-ohno1} 
  \\
  &&\sum_{n_1>\dots>n_{r}>0}
   \left\{
    \frac{q^{n_1a_1}\lambda}{[n_1]^{a_1}}
    -\frac{q^{n_1(a_1-1)}}{[n_1]^{a_1-1}}
   \right\}\frac{1}{[n_1]-q^{n_1}\lambda}
   \prod_{i=2}^s
    \frac{q^{(n_{c_{i}}-d_j)a_i}}
         {[n_{c_{i}}-d_j]^{a_i}}
   \prod_{j=2}^{r}\frac{1}{[n_{j}]-q^{n_{j}}\lambda}. 
   \nonumber
\end{eqnarray} 
Now we use the formula 
\begin{eqnarray*}
 \lefteqn{
  \left\{
    \frac{q^{na}\lambda}{[n]^{a}}
    -\frac{q^{n(a-1)}}{[n]^{a-1}}
   \right\}\frac{1}{[n]-q^{n}\lambda}
 }\\ 
 &&= 
  \left\{
    \frac{q^{(n-1)a}\lambda'}{[n-1]^{a}}
    -\frac{q^{(n-1)(a-1)}}{[n-1]^{a-1}}
   \right\}\frac{1}{[n]-q^{n}\lambda}+
  \frac{q^{(n-1)(a-1)}}{[n-1]^{a}}-\frac{q^{n(a-1)}}{[n]^{a}}. 
\end{eqnarray*}
Then we have 
\begin{eqnarray}
  \lefteqn{  \eqref{eq:pf-ohno1}=
   \lambda' \rho((a_1,0),b_1,\dots;\lambda)
   - \rho((a_1-1,0),b_1,\dots;\lambda)}
 \nonumber\\
   &&\qquad+
  \sum_{n_2>\dots>n_{r}>0}
  \prod_{i=2}^s
    \frac{q^{(n_{c_{i}}-d_j)a_i}}
         {[n_{c_{i}}-d_j]^{a_i}}
   \prod_{j=2}^{r}\frac{1}{[n_{j}]-q^{n_{j}}\lambda}
  \nonumber\\ 
   && \qquad \qquad \qquad  {}\times 
  \sum_{n_{1}=n_{2}+1}^{\infty} 
  \left\{ 
   \frac{q^{(n_{1}-1)(a_{1}-1)}}{[n_{1}-1]^{a_{1}}}
   -\frac{q^{n_{1}(a_{1}-1)}}{[n_{1}]^{a_{1}}}
  \right\}.
\label{eq:pf-ohno2} 
\end{eqnarray} 
{}From the equality 
\[
  \sum_{n_{1}=n_{2}+1}^{\infty} 
  \left\{ 
   \frac{q^{(n_{1}-1)(a_{1}-1)}}{[n_{1}-1]^{a_{1}}}
   -\frac{q^{n_{1}(a_{1}-1)}}{[n_{1}]^{a_{1}}}
  \right\}
  {}=
  \frac{q^{n_{2}(a_{1}-1)}}{[n_{2}]^{a_{1}}}=
  \frac{q^{n_{2}a_{1}}}{[n_{2}]^{a_{1}}}+(1-q)\frac{q^{n_{2}(a_{1}-1)}}{[n_{2}]^{a_{1}-1}}, 
\]
we obtain 
\begin{eqnarray*}
  \eqref{eq:pf-ohno2}&=&
   \lambda' \rho((a_1,0),b_1,\dots;\lambda)
   - \rho((a_1-1,0),b_1,\dots;\lambda) \\ 
 & &\quad {}+ 
   \rho((a_1,0),b_1-1,(a_2,0),\dots;\lambda) \\ 
 & &\qquad {}+
   (1-q)\rho((a_1-1,0),b_1-1,(a_2,0),\dots;\lambda). 
\end{eqnarray*} 
This completes the proof. 
\end{proof} 

\begin{proof}[Proof of Proposition~\ref{prop:difference}]
The left hand side of \eqref{eq:difference} is equal to 
\[
 \sum_{\eta \in I^{s}}
   (1-q)^{h(\eta)}(-\lambda)^{s-|\eta|+h(\eta)}
    \rho((a_1-\epsilon_1,0),b_1-\delta_1,
            \dots,(a_s-\epsilon_s,0),b_s-\delta_s;\lambda). 
\]
Take the sum over $(\epsilon_{1}, \delta_{1}) \in I$ 
by using Lemma \ref{lem:difference} \eqref{item:1a} and \eqref{item:1b}. 
Then we have 
\begin{eqnarray*}
\lefteqn{
   \sum_{\eta \in I^{s-1}}\sum_{\epsilon_{1}'=0, 1}
   (1-q)^{h(\eta)}(-\lambda)^{s-|\eta|-\epsilon_{1}'+h(\eta)}
   }\\ 
 &&{}\times 
    \rho((a_1-\epsilon_{1}',0),b_1,(a_{2}-\epsilon_{2}, 0), b_{2}-\delta_{2}, 
            \dots,(a_s-\epsilon_s,0),b_s-\delta_s;\lambda), 
\end{eqnarray*}
where $\eta=\{(\epsilon_{i}, \delta_{i})\}_{i=2}^{s}$. 
Next use \eqref{item:2} for $i=2, \ldots , s$ repeatedly, 
and apply \eqref{item:3a}. 
In the case of $b_{s}=1$, use \eqref{item:3b} at the last step. 
As a result we obtain 
\begin{eqnarray*}
 \lefteqn{
  \sum_{\eta' \in I^{s-1}}\sum_{\epsilon_{1}', \delta_{s+1}'=0, 1} 
   (1-q)^{h(\eta')}
    (-\lambda')^{s-|\eta'|-\epsilon_{1}'-\delta_{s+1}'+h(\eta')}
 }\\
 &&\times
     f(a_{1}-\epsilon_{1}',b_{1}-\delta_{2}',a_2-\epsilon_{2}',
      \dots,b_{s-1}-\delta_{s}',a_s-\epsilon_{s}',b_s-\delta_{s+1}';\lambda'), 
\end{eqnarray*} 
where $\eta'=\{(\delta_{i}', \epsilon_{i}')\}_{i=2}^{s}$. 
This is equal to the right hand side of \eqref{eq:difference}. 
\end{proof}

\begin{proof}[Proof of Theorem~\ref{thm:ohno}]
We prove \eqref{eq:f=g}
by induction on the value $\sum_{i}(a_{i}+b_{i})$. 
If $(a_1,b_1)=(1,1)$,
\[
 f(1,1;\lambda) = g(1,1;\lambda)
\]
is obvious from the definition.

Recall that 
\[
 g(a_{1}, b_{1}, \ldots , a_{s}, b_{s})=
 f(b_{s}, a_{s}, \ldots , b_{1}, a_{1}). 
\]
Hence the function $g$ also satisfies the equation \eqref{eq:difference}. 
Assume that \eqref{eq:f=g} holds
for the codes less than $(a_1,b_1,\dots,a_s,b_s)$.
{}From the induction hypothesis,
the difference of \eqref{eq:difference} for $f$ and $g$ gives 
\begin{eqnarray*}
\lefteqn{
 \lambda^s
 \bigl(
  f(a_1,b_1,\dots,a_s,b_s;\lambda) - g(a_1,b_1,\dots,a_s,b_s;\lambda)
 \bigr)}
 \\
 &&=
 {\lambda'}^s
 \bigl(
  f(a_1,b_1,\dots,a_s,b_s;\lambda') - g(a_1,b_1,\dots,a_s,b_s;\lambda')
 \bigr),
\end{eqnarray*}
i.e.\ the left hand side is invariant
under the variable transform $\lambda\mapsto q\lambda-1$.
On the other hand,
because of Proposition~\ref{prop:analytic},
we have 
\[
 \lambda^s
 \bigl(
  f(a_1,b_1,\dots,a_s,b_s;\lambda) - g(a_1,b_1,\dots,a_s,b_s;\lambda)
 \bigr)
 =
 \lambda^s\sum_{p=1}^\infty \frac{\tilde{C}_p}{[p]-q^p\lambda}
\]
for certain constants $\tilde{C}_{p}$. 
Using the invariance under $\lambda \mapsto q\lambda-1$, 
we find $\tilde{C}_p=0$ for all $p$.
\end{proof}

\subsection{The Ohno-Zagier Formula}

First we introduce the $q$-multiple polylogarithms: 
\begin{defn} 
For an index $\bfk =(k_{1}, \ldots , k_{r})$ 
the $q$-multiple polylogarithm (of one variable) $\Li_{\bfk}(t)$ 
is defined by 
\bea 
\Li_{\bfk}(t):=\sum_{n_{1}> \cdots >n_{r}>0}
\frac{t^{n_{1}}}{[n_{1}]^{k_{1}} \cdots [n_{r}]^{k_{r}}}. 
\label{def:q-multilog} 
\ena 
The right hand side of \eqref{def:q-multilog} is 
absolutely convergent if $|t|<1$. 
\end{defn} 

The $q$-multiple polylogarithms are related to $q$MZV as follows: 
\bea 
\Li_{k_{1}, \ldots , k_{r}}(q)&=&
\sum_{a_{1}=2}^{k_{1}} \sum_{a_{2}=1}^{k_{2}} \cdots \sum_{a_{r}=1}^{k_{r}} 
\binom{k_{1}-2}{a_{1}-2} 
\left\{ 
\prod_{j=2}^{r}\binom{k_{j}-1}{a_{j}-1} 
\right\} \nn \\ 
&& \qquad \qquad {}\times 
(1-q)^{\sum_{j=1}^{r}(k_{j}-a_{j})}  
\qzeta(a_{1}, \ldots , a_{r}). 
\label{eq:log-to-zeta} 
\ena 
Here $(k_{1}, \ldots , k_{r})$ is an admissible index. 

Now let us prove Theorem \ref{thm:ohno-zagier}. 
Denote by $I(k, r, s)$ 
the set of indices of weight $k$, depth $r$ and height $s$, 
and by $I_{0}(k, r, s)$ the subset consisting of admissible ones. 
Set 
\be 
G(k, r, s; t):=\sum_{\bfk \in I(k, r, s)}\Li_{\bfk}(t), \quad 
G_{0}(k, r, s; t):=\sum_{\bfk \in I_{0}(k, r, s)}\Li_{\bfk}(t). 
\en 
By definition we set $G(0, 0, 0; t)=1$, and 
$G(k, r, s; t)=0$ unless $k \ge r+s$ and $r \ge s \ge 0$. 
Consider the following generating functions 
\be 
&& 
\Phi:=\sum_{k, r, s \ge 0}G(k, r, s; t)u^{k-r-s}v^{r-s}w^{s}, \\ 
&& 
\Phi_{0}:=\sum_{k, r, s \ge 0}G_{0}(k, r, s; t)u^{k-r-s}v^{r-s}w^{s-1}. 
\en 
The function $\Phi_{0}$ is related to $q$MZV as follows: 
\begin{lem}\label{lem:phi-to-zeta} 
\be 
\Phi_{0}|_{t=q}= 
\frac{1}{1-(1-q)u}\sum_{k, r, s \ge 0}\left( 
\sum_{\bfk \in I_{0}(k, r, s)}\qzeta(\bfk) \right) 
x^{k-r-s}y^{r-s}z^{s-1}, 
\en 
where $x, y$ and $z$ are given by 
\bea 
x=\frac{u}{1-(1-q)u}, \quad y=\frac{v+(1-q)(w-uv)}{1-(1-q)u}, \quad 
z=\frac{w}{(1-(1-q)u)^{2}}. 
\label{eq:xyz-uvw}
\ena 
\end{lem} 
It is easy to prove this lemma from \eqref{eq:log-to-zeta}. 

To prove Theorem \ref{thm:ohno-zagier}
we derive a $q$-difference equation satisfied by $\Phi_{0}=\Phi_{0}(t)$. 
Denote by $\mathcal{D}_{q}$ the $q$-difference operator 
\bea 
(\mathcal{D}_{q}f)(t):=\frac{f(t)-f(qt)}{(1-q)t}. 
\label{def:q-diff} 
\ena 
Then the $q$-multiple polylogarithms satisfy 
\be 
\mathcal{D}_{q}\Li_{k_{1}, \ldots , k_{r}}(t)=\left\{ 
\begin{array}{ll} 
\displaystyle 
\frac{1}{t}\Li_{k_{1}-1, k_{2}, \ldots , k_{r}}(t), & k_{1} \ge 2, \\ 
\displaystyle
\frac{1}{1-t}\Li_{k_{2}, \ldots , k_{r}}(t), & k_{1}=1. 
\end{array} 
\right. 
\en 
These relations can be rewritten 
in terms of $G(k, r, s; t)$ and $G_{0}(k, r, s; t)$ as follows: 
\be 
&& 
\mathcal{D}_{q}G_{0}(k, r, s; t) \nn \\ 
&&\qquad {}= 
\frac{1}{t}\left( 
G(k-1, r, s-1; t)-G_{0}(k-1, r, s-1; t)+G_{0}(k-1, r, s;t) \right), \\ 
&& 
\mathcal{D}_{q}\left( 
G(k, r, s; t)-G_{0}(k, r, s; t) \right)= 
\frac{1}{1-t}G(k-1, r-1, s; t),  
\en 
or, in terms of generating functions,   
\be 
\mathcal{D}_{q}\Phi= 
\frac{1}{vt}\left( \Phi -1-w\Phi_{0} \right)+\frac{u}{t}\Phi_{0}, \quad 
\mathcal{D}_{q}\left( \Phi-w\Phi_{0} \right)= 
\frac{v}{1-t}\Phi. 
\en 
Eliminate $\Phi$ using the formula 
\be 
\mathcal{D}_{q}(tf(t))=qt \cdot \mathcal{D}_{q}f(t)+f(t). 
\en 
Then we obtain 
\bea 
qt(1-t)\mathcal{D}_{q}^{2}\Phi_{0}+ 
\left( (1-u)(1-t)-vt \right) \mathcal{D}_{q}\Phi_{0}+
(uv-w)\Phi_{0}=1. 
\label{eq:q-hyp} 
\ena 
This is an equation for the power series $\Phi_{0}=\Phi_{0}(t)$. 
Note that $\Phi_{0}(0)=0$. 
There is a unique solution to \eqref{eq:q-hyp} 
vanishing at $t=0$. 
To write down the solution we introduce the $q$-hypergeometric function 
\be 
\phi(a, b, c; t):=\sum_{n=0}^{\infty}t^{n}
\prod_{j=1}^{n}\frac{(1-aq^{j-1})(1-bq^{j-1})}{(1-q^{j})(1-cq^{j-1})} 
\en 
Then the solution is given by 
\bea 
\Phi_{0}(u, v, w; t)=\frac{1}{uv-w} \left( 
1-\phi(a, b, c; \frac{ct}{qab}) \right). 
\label{eq:solution} 
\ena 
Here $a, b$ and $c$ are defined in terms of $u, v$ and $w$ as follows: 
\bea
a=\frac{1}{1-(1-q)(u-\alpha_{0})}, \quad b=\frac{1}{1-(1-q)(u-\beta_{0})}, \quad 
c=\frac{q}{1-(1-q)u}, 
\label{eq:def-abc} 
\ena 
where $\alpha_{0}$ and $\beta_{0}$ are determined by 
\be 
\alpha_{0}+\beta_{0}=u+v, \quad \alpha_{0}\beta_{0}=w. 
\en 

Now we use Heine's $q$-analogue of Gauss' summation formula (see, e.g. \cite{MR91d:33034}):  
\be 
\phi(a, b, c; \frac{c}{ab})=
\prod_{n=0}^{\infty}
\frac{(1-a^{-1}cq^{n})(1-b^{-1}cq^{n})}
     {(1-cq^{n})(1-a^{-1}b^{-1}cq^{n})}. 
\en 
In our case, substituting \eqref{eq:def-abc}, we have 
\bea 
\phi(a, b, c; \frac{c}{ab})=\prod_{n=1}^{\infty} 
\frac{\displaystyle \left(1-\frac{q^{n}}{[n]}\alpha\right)
      \left(1-\frac{q^{n}}{[n]}\beta\right)}
     {\displaystyle \left(1-\frac{q^{n}}{[n]}x\right)
      \left(1-\frac{q^{n}}{[n]}y\right)}. 
\label{eq:product-form} 
\ena 
Here $x$ and $y$ are given by \eqref{eq:xyz-uvw}, and $\alpha$ and $\beta$ 
are defined by 
\be 
\alpha=\frac{\alpha_{0}}{1-(1-q)u}, \quad \beta=\frac{\beta_{0}}{1-(1-q)u}, 
\en 
or equivalently determined by 
\be 
\alpha+\beta=x+y+(q-1)(z-xy), \quad \alpha\beta=z. 
\en 
{}From Lemma \ref{lem:phi-to-zeta}, \eqref{eq:solution}, \eqref{eq:product-form} 
and the formula
\be 
\log \prod_{n=1}^{\infty}\left(1-\frac{q^{n}}{[n]}s \right)
&=& 
\frac{1}{q-1}\log( 1-s(q-1) ) \sum_{n=1}^{\infty}\frac{q^{n}}{[n]} \\ 
&&\qquad {}-
\sum_{n=2}^{\infty}\qzeta(n)\sum_{m=0}^{\infty}\frac{(q-1)^{m}}{m+n}s^{m+n}, 
\en 
we obtain Theorem \ref{thm:ohno-zagier}.

\section{Discussion}\label{sec:discussion} 
Now we consider the space 
spanned by the modified $q$MZV's \eqref{eq:modified-qMZV}. 
We regard the modified $q$MZV as a formal power series of $q$, and define 
subspaces of $\Q[[q]]$ by   
\[
Z_{k}:=\sum_{|\bfk|=k}\mathbb{Q}\qzetabar(\bfk),\quad
Z_{\le k}:=\sum_{2 \le |\bfk|\le k}\mathbb{Q}\qzetabar(\bfk). 
\]
Let us consider the problem to determine the dimension of the spaces over $\Q$. 
In principle, 
a lower bound for the dimension can be obtained as follows. 
Expand $q$MZV $\qzetabar(\bfk)$ as a power series of $q$: 
\[
\qzetabar(\bfk)=\sum_{n=0}^{\infty}a_{n}(\bfk)q^{n}, \quad 
a_{n}(\bfk) \in \mathbb{Z}_{\ge 0}. 
\]
Note that $a_{n}(\bfk)=0$ if $n < |\bfk|-1$. 
Recall that the number of admissible indices of weight $k$ is equal to $2^{k-2}$. 
Now consider the following square matrices 
\[
 A_{k}:=(a_{n}(\bfk))_{k-1 \le n \le k+2^{k-2}-2 \atop |\bfk|=k},
 \quad
 A_{\le k}:=(a_{n}(\bfk))_{1 \le n \le k+2^{k-2}-2 \atop
                                      2 \le |\bfk|\le k}. 
\]
Then we have 
\[
 \rank A_{k} \le \dim{Z_{k}},
 \quad
 \rank A_{\le k} \le \dim{Z_{\le k}}. 
\]
Let $d_k$ be the conjectured dimension
of the space of MZV's of weight $k$ over $\Q$ in \cite{MR96k:11110}.
In \cite{MR2003h:11073} it is shown that
$d_k$ gives the upper bound of the dimension of MZV's.

The values of $\rank A_{k}$ and $\rank A_{\le k}$ for $2 \le k \le 10$ are 
given as follows: 
\[
\begin{array}{@{\vrule width 1pt\ }c||r|r|r|r|r|r|r|r|r@{\ \vrule width 1pt}}
 \noalign{\hrule height 1pt}
 \hbox{weight}&
     2 & 3 & 4 & 5 & 6 & 7 & 8 &  9 & 10\\
 \noalign{\hrule height 1pt}
 d_k & 1 & 1 & 1 & 2 & 2 & 3 & 4 & 5 & 7\\
 \hline
 \hbox{$\rank A_{k}$}&
  1 & 1 & 2 & 3 & 6 & 9 & 18 & 29& 54\\
 \hline
 \hbox{By cyclic and Ohno}&
     1 & 1 & 2 & 3 & 6 & 9 & 18 & 30 & 56\\
 \noalign{\hrule height 1pt}
  \sum_{j\le k}d_j & 1 & 2 & 3 & 5 & 7 & 10 & 14 & 19 & 26\\
 \hline
 \hbox{$\rank A_{\le k}$}&
     1 & 2 & 4 & 7 & 11 & 18 & 27 & 42 & 63\\
 \hline
 \hbox{$\sum_{j\le k}\rank A_{j}$}&
    1 & 2 & 4 & 7 & 13 & 22 & 40 & 69 & 123\\
 \hline
  \noalign{\hrule height 1pt}
\end{array}
\]
Here the third row shows the upper bound for the dimension of $Z_{k}$ 
which follows from
the cyclic sum formula and the Ohno relation
\cite{Kaneko}.
We show some data of $a_{n}(\bfk)$ for $n\le 13$ and $|\bfk|\le 6$: 
\[
\begin{array}{@{\vrule width 1pt\ }c||rrrrrrrrrrrrr@{\ \vrule width 1pt}}
 \noalign{\hrule height 1pt}
\bfk & q^{1} & q^{2} & q^{3} & q^{4} & q^{5} & q^{6} & q^{7} & q^{8} & q^{9} & q^{10} & q^{11}&q^{12}&q^{13}\\
 \noalign{\hrule height 1pt}
(2) & 1& 3& 4& 7& 6& 12& 8& 15& 13& 18& 12&28&14\\
\hline 
(3) & 0& 1& 3& 7& 10& 19& 21& 35& 39& 56& 55&91&78\\
\hline 
(4) & 0&0& 1& 4& 10& 21& 35& 60& 85& 130& 165&245&286\\
(3,1) & 0&0&0& 1& 1& 6& 5& 15& 18& 31& 30&70&55\\
\hline 
(5) & 0&0&0& 1& 5& 15& 35& 71& 126& 215& 330&511&715\\
(4,1) & 0&0&0&0&0& 1& 1& 5& 7& 16& 17&47&42\\
(3,2) & 0&0&0&0& 1& 2& 7& 13& 24& 42& 69&97&149\\
\hline 
(6) & 0&0&0&0& 1& 6& 21& 56& 126& 253& 462&798&1287\\
(5,1) & 0&0&0&0&0&0&0& 1& 1& 6& 6&23&22\\
(4,2) & 0&0&0&0&0&0& 1& 2& 7& 13& 30&45&88\\
(3,3) & 0&0&0&0&0&1&3&10&22&47&85&154&244\\
(4,1,1) &0&0&0&0&0&0&0&0& 1& 1& 2 & 9 & 9\\
(3,2,1) & 0&0&0&0&0&0&0& 1& 1& 4& 9 &14 &23\\
 \noalign{\hrule height 1pt}
\end{array}
\]
Any $q$MZV of weight less than $7$ is given by a linear combination of 
ones in the list above from the cyclic sum formula and the Ohno relation.

Now we discuss two points. 
First, at weight $9$ there is a gap between $\rank A_{9}$ and 
the upper bound.  
This shows a possibility that our linear relations are not sufficient to determine 
the dimension of $Z_{k}$. 
In fact, the following equality holds up to $q^{269}$: 
\begin{eqnarray*}
4 \qzeta(7,2) +6 \qzeta(6,3) - \qzeta(5,4) - \qzeta(4,5) -6 \qzeta(6,2,1) -6 \qzeta(6,1,2)&&\\
 -2 \qzeta(5,3,1) -7 \qzeta(5,2,2)-3 \qzeta(5,1,3) +2 \qzeta(4,4,1) - \qzeta(4,3,2)&&\\
 +\qzeta(3,5,1) +\qzeta(3,2,4) -3 \qzeta(2,5,2) +2 \qzeta(5,2,1,1)&&\\
+2 \qzeta(5,1,2,1) +2 \qzeta(5,1,1,2) + \qzeta(3,3,1,2) - \qzeta(3,2,3,1)&&\\
 -4 \qzeta(3,2,2,2) - \qzeta(3,2,1,3) -2 \qzeta(2,2,3,2) +\qzeta(2,1,3,3) &\stackrel{?}{=}&0
\end{eqnarray*}
It is linearly independent of 
the cyclic sum formula and the Ohno relation,
and is
correct by substituting $\qzeta$ with $\zeta$
\cite{Kaneko}.

Next, if the weight is greater than $5$, the equality 
$\rank A_{\le k}=\sum_{j=2}^{k}\rank A_{j}$ breaks. 
This shows that there may exist linear relations 
among the modified $q$MZV's of different weights. 
For example,
as observed in the data above, 
the following relations seem to hold:
\begin{eqnarray}
 && -\qzetabar(3,1) + \qzetabar(5) -3\qzetabar(4,1)
 -3\qzetabar(3,2)
 \nonumber\\
 &&\hspace{5em}
  - \qzetabar(6) -3 \qzetabar(4,2)
 +6\qzetabar(3,3) \stackrel{?}{=} 0,
 \nonumber\\
 &&
  -2\qzetabar(3,1) + 2\qzetabar(5) -6 \qzetabar(4,1) -9 \qzetabar(3,2)
  \label{eq:weight6}\\
 &&\hspace{5em}
 + \qzetabar(6) -12 \qzetabar(4,2) - 3 \qzetabar(4,1,1)+
 3\qzetabar(3,2,1)
 \stackrel{?}{=} 0.
 \nonumber
\end{eqnarray}
Unfortunately, they are not derived from the relations proved in this paper. 
An interesting point is that 
we obtain some relations of MZV's from the formulae above. 
Multiply $(1-q)^6$ and take the limit $q \to 1$ in \eqref{eq:weight6}. 
Then the terms of weight less than $6$ vanish and as a result we find 
\begin{eqnarray*}
 - \zeta(6) -3 \zeta(4,2)
 +6\zeta(3,3) &=& 0,\\
 \zeta(6) -12 \zeta(4,2) - 3 \zeta(4,1,1)+
 3\zeta(3,2,1)
 &=& 0.
\end{eqnarray*}
These relations for MZV's are correct 
\cite{Kaneko}.

\section*{Aknowledgements}
The authors wish to thank Kaneko Masanobu and Ohno Yasuo
for many helpful informations and comments.

\def\cprime{$'$} \def\cprime{$'$}

\begin{tabular}{l}
 OKUDA Jun-ichi,\\
 Department of Mathematical Sciences\\
 Graduate School of Science and Engineering\\
 Waseda University,
 Tokyo 169-8555, Japan\\
 okuda@gm.math.waseda.ac.jp\\
 \\
 TAKEYAMA Yoshihiro,\\
 Department of Mathematics,\\
 Graduate School of Science,\\
 Kyoto University, Kyoto 606-8502, Japan\\
 takeyama@math.kyoto-u.ac.jp
\end{tabular}
\end{document}